\documentclass[a4paper,12pt]{article}
\usepackage{amsmath}
\usepackage{amsthm}
\usepackage{amssymb, amscd}
\usepackage{epic, eepic}
\usepackage[all]{xy}
\usepackage{mathrsfs}
\usepackage{enumerate}
\title{Complete Intersections of Two Quadrics and Galois Cohomology}
\author{Yasuhiro Ishitsuka}
\date{\today}

\newcommand{\Q}{\mathbb{Q}}

\newcommand{\Z}{\mathbb{Z}}

\newcommand{\Oh}{\mathcal{O}}

\newcommand{\p}[1]{\mathbb{P}^{#1}}

\def\dim{\mathop{\mathrm{dim}}\nolimits}

\def\id{\mathop{\mathrm{id}}\nolimits}

\def\Tr{\mathop{\mathrm{Tr}}\nolimits}
\def\det{\mathop{\mathrm{det}}\nolimits}
\def\N{\mathop{\mathrm{N}}\nolimits}

\def\Ker{\mathop{\mathrm{ker}}\nolimits}

\def\Hom{\mathop{\mathrm{Hom}}\nolimits}

\def\End{\mathop{\mathrm{End}}\nolimits}
\def\Aut{\mathop{\mathrm{Aut}}\nolimits}
\def\GL{\mathop{\mathrm{GL}}\nolimits}

\def\Gal{\mathop{\mathrm{Gal}}\nolimits}
\def\Res{\mathop{\mathrm{Res}}\nolimits}
\makeatletter
\def\iddots{\mathinner{\mkern1mu\raise\p@
    \hbox{.}\mkern2mu\raise4\p@\hbox{.}\mkern2mu
    \raise7\p@\vbox{\kern7\p@\hbox{.}}\mkern1mu}}
\def\adots{\mathinner{\mkern2mu\raise\p@\hbox{.} 
 \mkern2mu\raise4\p@\hbox{.}\mkern1mu
 \raise7\p@\vbox{\kern7\p@\hbox{.}}\mkern1mu}}
\makeatother

\theoremstyle{plain}
\newtheorem{thm}{Theorem}[section]
\newtheorem{prop}[thm]{Proposition}
\newtheorem{cor}[thm]{Corollary}
\newtheorem{lem}[thm]{Lemma}

\theoremstyle{definition}
\newtheorem{dfn}[thm]{Definition}
\newtheorem*{prf}{Proof}

\newtheorem{rmk}[thm]{Remark}
\newtheorem{exa}[thm]{Example}
\newcommand{\sq}{\qquad $\square$}

\newcommand{\mc}[1]{\mathcal{#1}}

\newcommand{\ms}[1]{\mathsf{#1}}

\newcommand{\isoc}{\overset{\sim}{\rightarrow}}

\newcommand{\Sym}[2]{{\mathrm{Sym}_{#1} #2}} 
\newcommand{\sym}[2]{{\mathrm{Sym}^{#1} #2}} 

\def\Spec{\mathop{\mathrm{Spec}}\nolimits}

\begin{document}
\maketitle
\begin{abstract}
For each nonsingular hyperelliptic curve of arbitrary genus with a rational Weierstrass point,
we construct explicitly a natural map from the Galois cohomology of 2-torsion points of the Jacobian variety
of the curve to the set of projective equivalence classes of nonsingular complete intersections of two quadrics. 
We also investigate a necessary condition for an element in the Galois cohomology
 to be in the image of a rational point of the Jacobian variety.
This gives a generalization of the results of Flynn and Skorobogatov.
\end{abstract}

\section*{Introduction}
Flynn and Skorobogatov constructed maps from the Galois cohomology of
 2-torsion points of Jacobian varieties of nonsingular hyperelliptic curves of genus 2
 with rational Weierstrass points to the set of isomorphism classes 
of del Pezzo surfaces of degree 4 (see \cite{Flynn}, \cite{Skor2}).
They also proved that the images of the maps for all such curves 
cover the whole set of isomorphism classes of del Pezzo surfaces of degree 4,
and characterized the images of the trivial elements. 
In this paper, we generalize these results to nonsingular hyperelliptic curves of arbitrary genus with rational Weierstrass points.
We use the methods of Wood (\cite[\S 5]{Wood}).


	Let $k$ be a field of characteristic not equal to 2 with separable closure $k^s$,
	and $G_k := \Gal(k^s/k)$ the absolute Galois group.
	A separable polynomial $f(t)$ of odd degree $n = 2m+1 \geq 3$ with coefficients in $k$
	defines an \'etale $k$-algebra $L := k[t]/\left( f(t) \right)$ of dimension $2m+1,$
	and a nonsingular projective hyperelliptic curve $C$ over $k$ of genus $m$ whose affine model is given by the equation \[y^2 = f(t).\] 
	Let $J_C$  be the Jacobian variety of $C$. The group of 2-torsion points $J_C[2](k^s)$ is isomorphic to
	the group of $k^s$-valued points of a group scheme 
\[
G := \Res_{L/k}(\mu_2) / \Delta(\mu_2)
\]
	as $G_k$-modules, where $\Res_{L/k}$ denotes Weil restriction and
\[
 \Delta \colon \mu_2 \hookrightarrow \Res_{L/k}(\mu_2) 
\] 
	denotes the diagonal embedding.
	Since the short exact sequence
	\[\begin{CD}
		0 @>>> \mu_2 @>\Delta>> \Res_{L/k}(\mu_2) @>>> G \cong J_C[2] @>>> 0
	\end{CD}\]
	is split, we have isomorphisms \[H^1 \left( G_k, J_C[2] \right) \cong H^1 \left( G_k, G \right) \cong L^\times / k^\times L^{\times 2}.\]
	
	The main goal of this paper is the following theorem:
\begin{thm}[See Proposition \ref{main1} and Theorem \ref{quasi}]\label{main}	
	Let $m \geq 1$ be a positive integer, and $k$ a field of characteristic not equal to $2.$ 
	We put $n = 2m+1,$ and assume $\# k \geq n$.\\
	(i) Let $C$ be a nonsingular projective hyperelliptic curve of genus $m$ defined over $k$
	whose affine model is given by $y^2= f(t)$ where $f(t) \in k[t]$ is a separable polynomial of degree $n = 2m+1.$
	Let $\ms{ciq}_n$ be the set of projective equivalence classes of nonsingular polarized complete intersections of two quadrics in $\p{n-1}_k$ over $k.$
	Then, there exists a natural map 
\[
 i_C \colon H^1 \left( G_k, J_C[2] \right) \to 
\ms{ciq}_n
\]
	such that $i_C(0)$ is a variety containing a linear $k$-subvariety isomorphic to $\p{m-1}_k$ (we call it quasi-split).\\
	(ii)
	Conversely, for a quasi-split variety $X$, there exists a nonsingular projective hyperelliptic curve $C$ of genus $m$
	defined over $k$ with a rational Weierstrass point such that $i_C(0) = [X].$ \\
	(iii)
	Moreover, for any projective equivalence class $[X]$ of nonsingular polarized complete intersection $X$
	of two quadrics in $\p{n-1}_k$ over $k$, 
	there exists a nonsingular projective hyperelliptic curve $C$
	of genus $m$ defined over $k$ with a rational Weierstrass point
	and $\eta \in H^1 \left( G_k, J_C[2] \right)$
	such that $i_C(\eta) =[X].$
\end{thm}
\begin{rmk}
	We also determine the fiber $i_C^{-1}([X])$ in Proposition \ref{Ru}.
	When $n \geq 5,$ the projective equivalence classes coincide with isomorphism classes of those varieties (see Remark \ref{isom}).
\end{rmk}
	Recall that we have the following short exact sequence
	\[ \begin{CD}
		0 @>>> J_C(k)/2J_C(k) @>{\delta}>> H^1 \left( G_k, J_C[2] \right) @>>> H^1 \left( G_k, J_C \right)[2] @>>> 0
	\end{CD}.\]
	It is a natural and interesting question to ask which elements in $\ms{ciq}_n$ come from $k$-rational points of the Jacobian variety $J_C.$
 	On the other hand, the map $i_C$ is not injective in general. Moreover, it is easy to find two curves $C, C'$ over $k$
 	which are not $\overline{k}$-isomorphic but define the same class 
	$[X] \in \ms{ciq}_n$ such that $[X]$ is in the image of $i_C \circ \delta$,
	but not in the image of $i_{C'} \circ \delta.$
 	Therefore it seems difficult to determine the image of $i_C \circ \delta$ exactly.
 	The following theorem gives a necessary condition for an element in $\ms{ciq}_n$ to come 
	from a $k$-rational point of $J_C.$

  
 \begin{thm}[see Theorem \ref{Jp}]\label{main2}
	For a $k$-rational point $P \in J_C(k),$ we put $i_C(\delta(P)) = [X].$
	Then, the variety $X$ contains a linear $k$-subvariety isomorphic to $\displaystyle \p{\lfloor \frac{m-1}{2} \rfloor}_k.$
 \end{thm}

Note that Theorem \ref{main2} is not a sufficirnt condition.
When $n=5,$ Theorem \ref{main} and Theorem \ref{main2} are proved by Flynn and Skorobogatov (\cite{Flynn}, \cite{Skor2}).
Flynn used his result to construct nontrivial elements for the Shafarevich--Tate group of the Jacobian of hyperelliptic curve of genus 2.
We give a different proof of their results when $n=5.$
We use the method of Wood \cite[\S 5]{Wood} relating elements in $\ms{ciq}_n$ and monogenic \'etale algebras.

In Section 1, we define $\ms{ciq}_n$ and other sets $\ms{orb}_n, \ms{gc}_n,$
 and $\ms{hec}_n.$ Then we construct explicit bijections between them. 
In Section 2, we define the map $i_C$ using the bijections in Section 1,
 and show Theorems \ref{main}, \ref{main2}. 
In Section 3, we give examples of low genus and compare our results with Skorobogatov's results.
\subsection*{Notation}
	In this paper, we fix a field $k$ of characteristic not equal to 2 and its algebraic closure $\overline{k}$.
	Denote the separable closure of $k$ in $\overline{k}$ by $k^s,$ 
	and the absolute Galois group by $G_k = \Gal (k^s/k).$ 
	For $\theta \in k^s$ and $A=\begin{pmatrix} a & b \\ c& d \end{pmatrix} \in \GL_2(k),$
	we denote $\displaystyle \frac{a\theta + b}{c\theta + d} \in k^s$ by $A \cdot \theta.$
	
	Let $m \geq 1$ be a positive integer and $n = 2m+1.$ We assume $\# k \geq n.$ 
	We need this assumption only to construct a map $\ms{orb}_n \to \ms{gc}_n$ in Section 1.2. 
	We use the following notation: $V$ an $n$-dimensional $k$-vector space, 
	$\{v_i\}_{0 \leq i \leq n-1}$ a basis of $V,$ $W$ a two-dimensional $k$-vector space,
	$\{w_i\}_{i=0, 1}$ a basis of $W.$ 
	Also we use $V^*$ (resp. $W^*$) for the dual space of $V$ (resp. $W$), 
	and $\{v_i^*\}$ (resp. $\{w_i^*\}$) for the dual basis with respect to $\{v_i\}$ (resp. $\{w_i\}$). 
	For $r \geq 2$, we denote by $\Sym r V$ the $r$-ic symmetric subspace of $\bigotimes^r V$, 
	and $\sym r V$ the $r$-ic symmetric quotient of $\bigotimes^r V$. 
	We use similar notation for the alternating subspace $\wedge_r V$ 
	and the alternating quotient space $\wedge^r V.$ 
	When we say $L$ is a $k$-algebra of degree $n,$ 
	we mean $L$ is a commutative $k$-algebra whose dimension is $n$ as a $k$-vector space.
	Two subvarieties $X, Y$ of $\p{n}_k$ are said to be projectively equivalent if there is an automorphism $\p{n}_k \isoc \p{n}_k$
	which induces an isomorphism $X \isoc Y.$
 \subsection*{Acknowledgements}
	The author is very grateful to Professor Tetsushi Ito for much advice, useful remarks, and warmful encouragement during this work.
 
\section{Definition of $\ms{ciq}_n, \ms{orb}_n, \ms{gc}_n,$ and $\ms{hec}_n$}
In this section, we will define four sets, $\ms{ciq}_n, \ms{orb}_n, \ms{gc}_n,$ and $\ms{hec}_n,$
 and construct bijections between them.
To describe the main results, we only need to define two sets, the projective equivalence classes of 
polarized complete intersections of two quadrics $\ms{ciq}_n$ 
and a set related to hyperelliptic curves and Galois cohomology. But we use mainly the third set $\ms{gc}_n$,
which is defined in terms of \'etale $k$-algebras, and relatively easy to describe and compute.

\subsection{Definition and relation of geometric objects $\ms{ciq}_n$ and orbital objects $\ms{orb}_n$}
Before discussion, we recall some definitions from \cite{Reid}. 
Let $V$ be an $n$-dimensional $k$-vector space. 
Fix a two-dimensional subspace $W$ of the space of quadratic forms $\Sym2V^*$, 
which we call a {\it pencil of quadratic forms}. 
Each quadratic form defines a quadric in $\p{}(V)$, so $W$ defines a {\it pencil of quadrics}. 
We identify them and abuse the notation.

For a quadratic form $Q \in \Sym2 V^*,$ 
we define the {\it degeneracy} as $\dim_k \Ker (Q: V \rightarrow V^*).$
A quadratic form $Q \in \Sym2V^*$ is {\it nondegenerate} if the degeneracy of $Q$ is zero.
 Two pencils $W \subset \Sym2 V^*, W' \subset \Sym2 V'^*$ are {\it projectively equivalent} 
if there exists a $k$-linear isomorphism $V \isoc V'$ which induces an isomorphism $W \isoc W'.$
\begin{dfn}\label{fir}
	(a) A pencil of quadratic forms $W \subset \Sym2 V^*$ is called {\it nonsingular} 
if it satisfies the following two conditions: 
\\(i) the degeneracy of a quadratic form $Q \in W \otimes_k \overline{k}$ is at most one, 
\\(ii) $\bigcap_{Q \in W \otimes_k \overline{k}} \Ker(Q) = \{0\}.$ \\
	(b) The {\it determinant form} $\det \in \sym{n} W^* \otimes \left( \wedge^n V^* \right) ^{\otimes 2}$
 of a pencil of quadratic forms $W \subset \Sym2 V^*$ is 
$$W \ni w \mapsto \wedge^n w \in \left( \wedge^n V^* \right) ^{\otimes 2}.$$
\end{dfn}
The following theorem from \cite{Reid} says that two quadrics 
which intersect completely and smoothly can be identified 
with a basis of a nonsingular pencil of quadratic forms.
	\begin{prop}[{\cite[Proposition 2.1]{Reid}}] \label{Reid}
		Let $V$ be an $n$-dimensional $k$-vector space, 
		$W \subset \Sym2 V^*$ a pencil of quadrics in $\p{}(V)$, 
		and $X := \bigcap_{Q \in W \otimes_k \overline{k}} Q \subset \p{}(V \otimes_k \overline{k}).$
		Then the following conditions are equivalent:
		\begin{description}
		\setlength{\parskip}{0pt}
		\setlength{\itemsep}{0.0pt}
			\item[(a)] $X$ is nonsingular and codimension two in $\p{}(V)$
							(this means $X$ is a nonsingular complete intersection of two quadrics over $k$).
			\item[(b)] The pencil $W$ is nonsingular in the sense of Definition \ref{fir}.
			\item[(c)] The determinant form $\det \in \sym{n} W^* \otimes \left( \wedge^n V^* \right) ^{\otimes 2}$
							has nonzero discriminant.
		\end{description}
		If one considers over $k = k^s$, these conditions are equivalent to another condition:
		for any basis $\{ w_0, w_1 \}$ of $W$ such that $w_0$ is a nondegenerate quadratic form on $V$,
		there exists a basis $\{ v_i \}$ of $V$ and $\lambda_i \in k^s \; (0 \leq i \leq n-1)$ 
		such that for any $x_i \in k^s \; (0 \leq i \leq n-1),$ 
			\begin{align*}
				w_0 \left( \sum_{i=0}^{n-1} x_iv_i \right) &= \sum_{i=0}^{n-1} x_i^2 \\
				w_1 \left( \sum_{i=0}^{n-1} x_iv_i \right) &= \sum_{i=0}^{n-1} \lambda_ix_i^2,
			\end{align*}
			and $\lambda_i \neq \lambda_j$ for $i \neq j$. Furthermore, the basis $\{v_i\}$ is unique 
			up to permutation and change of signs in front of the $v_i$.
	\end{prop}
\begin{dfn}
 We define $\ms{ciq}_n$ as the set of projective equivalence classes of 
nonsingular polarized complete intersections of two quadrics over $k$ in $\p{}(V).$
Here, the polarization of $X$ means the embedding $X \hookrightarrow \p{}(V).$
\end{dfn}
\begin{rmk}\label{isom}
	Let $X$ be a nonsingular complete intersection of two quadrics over $k$ in $\p{n-1}_k.$
	If $n = 5,$ $X$ is a del Pezzo surface of degree four, and the embedding $X \hookrightarrow \p4_k$ is the canonical embedding.
	Hence, two del Pezzo surfaces of degree four $X, X'$ are projective equivalent if and only if $X, X'$ are $k$-isomorphic.
	If $n \geq 7,$ the same conclusion is obtained from the fact that $\mathrm{Pic}(X) =\Z.$ 
	So if $n \geq 5,$ the set $\ms{ciq}_n$ coincides with the set of isomorphism classes of nonsingular complete intersections of two quadrics.
	
	But when $n=3,$ this is not true: the elements of $\ms{ciq}_3$ are four points subschemes of $\p2_k.$
\end{rmk}
 Let $X$ be a nonsingular complete intersection of two quadrics in $\p{}(V).$
Then the subspace $W$ of $\Sym2 V^*$ defined by 
\begin{equation*}
\begin{split}
	W &:= \Ker \left( H^0 \left( \p{}(V), \Oh_{\p{}(V)}(2) \right) \rightarrow H^0 \left( X, \Oh_X(2) \right) \right) \\
	&\subset \Sym2 H^0 \left( \p{}(V), \Oh_{\p{}(V)}(1) \right) = \Sym2 V^*
\end{split}
\end{equation*}
　is two-dimensional and defines a nonsingular pencil of quadratic forms.
	If $X, X'$ are projectively equivalent, two pencils defined by $X, X'$ are projectively equivalent.



 Conversely, a nonsingular pencil of quadratic forms $W \subset \Sym2 V^*$
 defines a nonsingular complete intersection of two quadrics $X$ by Proposition \ref{Reid}. 
Two pencils $W, W'$ are projectively equivalent
 if and only if a coordinate change of $\p{}(V)$ induces an isomorphism $W \isoc W'.$ 
Hence all elements in a projective equivalence class of $W$ defines varieties isomorphic to $X$ 
and the polarization $X \hookrightarrow \p{}(V).$ These two maps are inverses of each other.

Furthermore, the set of projective equivalence classes of nonsingular pencils of quadratic forms
 can be identified with another set. 
Let us consider the space of pairs of quadratic forms 
 \[
k^2 \otimes \Sym2 k^n.
\]

 This space has a natural $\GL_2(k) \times \GL_n(k)$-action. 
Note that if $A, B \in \Sym2 k^n$ are linearly independent, 
the $\GL_2(k)$-orbit of $(A, B) \in k^2 \otimes \Sym2 k^n$ can be identified with
the two-dimensional subspace $W \subset \Sym2 k^n$ generated by $A$ and $B,$
and the $\GL_n(k)$-orbit of $W$ can be identified with the projective equivalence class of $W.$
Any pencils of quadratic forms are projectively equivalent to a pencil in $\Sym2 k^n.$

The {\it determinant form} of $(A, B) \in k^2 \otimes \Sym2 k^n$ is 
\[
 \det (Ax-By) \in \sym{n}k^2.
\]
We say a $ \GL_2(k) \times \GL_n(k)$-orbit of $(A, B) \in k^2 \otimes \Sym2 k^n$ is {\it stable}
when its determinant form has nonzero discriminant. 

\begin{dfn}
 We define $\ms{orb}_n$ as the set of stable  $\GL_2(k) \times \GL_n(k)$-orbits of $k^2 \otimes \Sym2 k^n.$
\end{dfn}

The stability is equivalent to the non-singularity of the pencil of quadratic forms $W \subset \Sym2 k^n.$
So the set of projective equivalence classes of nonsingular pencils of quadratic forms can be identified with $\ms{orb}_n.$

 Consequently, we construct the following bijections:
 \begin{equation}\label{first}
\ms{ciq}_n \overset{\sim}{\leftrightarrow} \left\{ \mbox{
\begin{tabular}{l}
projective equivalence classes\\ of nonsingular pencils of\\ quadrics in $\p{n-1}_k$ over $k$
\end{tabular} } \right\} \overset{\sim}{\leftrightarrow} \ms{orb}_n.
\end{equation}
 
 Next, we discuss the characteristic schemes.
 Recall that closed subschemes $S \subset \p{}(W)$ and $S' \subset \p{}(W')$ are 
{\it projectively equivalent} if there exists a $k$-isomorphism $W \isoc W'$
 which induces isomorphisms $\p{}(W) \isoc \p{}(W')$ and $S \isoc S'.$
\begin{dfn}
	The {\it characteristic scheme} of a $\GL_2(k) \times \GL_n(k)$-orbit of $(A, B) \in k^s \otimes \Sym2 k^n$
	is the projective equivalence class of the subscheme of $\p{1}$ where the determinant form vanishes.
	The {\it characteristic scheme} of a pencil $W$ is the projective equivalence class of the subscheme of $\p{}(W)$ 
	where the determinant form $\det \in \sym{n}W^* \otimes \left( \wedge^n V^* \right) ^{\otimes 2}$ vanishes.
\end{dfn}
	Strictly speaking, both characteristic schemes are not schemes, but projective equivalence classes of subschemes of $\p1_k.$
	It is easy to see that the bijections (\ref{first}) preserve  the characteristic schemes.

\subsection{Definition of algebraic objects $\ms{gc}_n$ and their relation to $\ms{orb}_n$}
	In this subsection, we define the third set $\ms{gc}_n$ and discuss its relation to the set $\ms{orb}_n.$
	The definition of $\ms{gc}_n$ is inspired by the results of Wood (\cite[\S 5]{Wood}).
	\begin{dfn}
		Let us consider triples $(L, \theta, \alpha)$ satisfying the following properties:
	\begin{itemize}
		\setlength{\parskip}{0pt}
		\setlength{\itemsep}{0.0pt}
		\item $L$ is an \'etale $k$-algebra of degree $n$,
		\item $\theta \in L$ is a generator of $L$ as a $k$-algebra,
		\item $\alpha $ is an element of $L^\times / k^\times L^{\times 2}$.
	\end{itemize}
		Two triples $(L, \theta, \alpha), (L', \theta', \alpha')$ are said to be {\it equivalent} if there exists a $k$-algebra isomorphism $L \isoc L'$ and 
$A= \left( \begin{matrix}
a&b\\
c&d
\end{matrix} \right) \in \GL_2(k)$
 such that $A \cdot \theta = \displaystyle \frac{a \theta+b}{c \theta+d} \mapsto \theta'$ and $\alpha \mapsto \alpha'.$
 We define $\ms{gc}_n$ as the set of equivalence classes of triples $(L, \theta, \alpha)$ as above.
\end{dfn}
  This equivalence relation can be separated into a combination of two equivalence relations:
 M\"{o}bius transformation of the generator $\theta \mapsto \displaystyle \frac{a \theta+b}{c \theta+d}$ 
and an isomorphism of $k$-algebras $L \isoc L'$. 
A M\"{o}bius transformation changes the representation of $\alpha$ 
as a polynomial of the generator $\theta$, while an isomorphism of $k$-algebras does not. 
 We note that 
\[
H^1 \left( G_k, \Res_{L/k}(\mu_2)/\Delta(\mu_2) \right) \cong L^\times / k^\times L^{\times 2}
\]
 because $n$ is odd.
 
 \vspace{5pt}
 Now we shall construct a bijection between $\ms{orb}_n$ and $\ms{gc}_n$. 
First, we construct a map $\ms{orb}_n \rightarrow \ms{gc}_n$. 
In fact, this step is the hardest part of the construction of $i_C.$
 Take an element $W \in \ms{orb}_n$ and fix a representative $(w_0, w_1) \in k^2 \otimes \Sym2 k^n$ of $W$ 
such that $w_0$ is a nondegenerate quadratic form on $V=k^n$ 
(in this choice, we use the assumption $\# k \geq n$). 
Consider $W$ as a two-dimensional $k$-linear subspace of $\Hom_k \left( V, V^* \right) ( \cong V^* \otimes V^*)$ and put
\begin{align*}
 w_0^{-1}W := \{w_0^{-1} \circ w \, | \, w \in W \} \subset \End V.
\end{align*}
The space $w_0^{-1}W \subset \End V$ is a two-dimensional $k$-vector space, and 
\[
	\{\id = w_0^{-1} \circ w_0, \quad \theta := w_0^{-1} \circ w_1 \}
\]
is a basis of $w_0^{-1}W.$ We write $L_W$ as a $k$-subalgebra of $\End V$ generated by $\theta$ (or, equivalently, $w_0^{-1}W$).
	\begin{lem}\label{fr}
		The $k$-algebra $L_W$ is an \'etale $k$-algebra of degree $n.$ Moreover, $V$ is a free $L_W$-module of rank one.
	\end{lem}
	\begin{prf}
		Write the characteristic polynomial of $\theta \in \End V$ as $P_\theta(t):= \det \left( t-\theta ; V \right)$. 
		By definition, \[P_\theta(t) = \det(t - w_0^{-1} \circ w_1 ; V).\]
		Since $\det \in \sym{n} W^* \otimes \left( \wedge^n V^* \right)^{\otimes 2}$
		has nonzero discriminant by Proposition \ref{Reid},
		$P_\theta(t)$ is a separable polynomial and
		the minimal polynomial of $\theta$ coincides with $P_\theta(t).$
		Hence, we see $L_W \cong k[t] / (P_\theta(t))$ and it is an \'etale $k$-algebra of degree $n.$
		
		Then $L_W$ can be written as a product of fields $\prod_{s=1}^{r} K_s$, where $K_s$ is a finite separable extension of $k$.
		We can take idempotent elements $e_s \in L_W$ such that 
		\[ L_W e_s = K_s, \quad e_se_t=0 \quad (s \neq t), \quad \mbox{and} \quad \sum_{s=1}^r e_s = 1_{L_W}.\]
		For each $s,$ we take vectors $v'_s \in V$ such that $e_s v'_s \neq 0.$ We put $v' := \sum_{s=1}^{r} e_s v'_s .$ 
		
		If an element $\ell \in L_W$ satisfies $\ell v' = 0,$ we find $\ell e_s v' = (\ell e_s) (e_sv'_s) =0$ for all $s.$
		Since $e_sv'_s \neq 0$ and $K_s$ is a field, the element $\ell e_s \in K_s \subset L_W$ is zero for all $s.$
		So $\ell = \sum_{s=1}^r \ell e_s = 0$ in $L_W.$ 
		By comparing dimensions as $k$-vector spaces, we conclude that $v'$ defines a $k$-isomorphism 
\begin{align*}
	L_W &\isoc V \quad \\ \quad \ell &\mapsto \ell v'.
\end{align*}
		Hence $V$ is a free $L_W$-module of rank one. \hfill \sq
	\end{prf}
	From Lemma \ref{fr}, we fix an $L_W$-module isomorphism $\gamma \colon L_W \isoc V.$
	Then, for each $w \in W$, we define $\widetilde{w} := (\gamma^* \otimes \gamma^*) (w) \in \Sym2 L_W^*.$
	Put $\widetilde{W}$ for the two-dimensional $k$-subspace of $\Sym2 L_W^*$ consisting of such elements. 
	
	Recall that the characteristic of $k$ is different from 2. 
	We identify quadratic forms on a $k$-vector space and symmetric bilinear forms on it.
	In particular, each element of $\widetilde{W}$ (resp. $W$) defines a symmetric bilinear form on $L_W$ (resp. $V$).  
\begin{lem} \label{de}
All symmetric $k$-bilinear forms $\widetilde{w} \in \widetilde{W} \subset \Sym2 L_W^*$ satisfy
	\begin{equation}\label{wdt}
		\widetilde{w}(\theta x , y) = \widetilde{w}(x, \theta y) \quad (\mbox{for } x, y \in L_W).
	\end{equation}
\end{lem}
\begin{prf}
	Since $\gamma$ is an isomorphism of $L_W$-modules, $\gamma(\theta x) =  \theta \gamma(x).$
	By regarding $\Sym2 V^* \subset V^* \otimes V^*$ as the $k$-subspace of $\Hom_k \left( V, V^* \right),$ we calculate
	\begin{align*}
		\widetilde{w}(\theta x, y) &= w \left( \theta\gamma(x), \gamma(y) \right)\\
		&= w \left( \theta \gamma(x) \right) \left( \gamma(y) \right).
	\end{align*}
	Because $\theta = w_0^{-1} \circ w_1$ and we can write $w=bw_0+aw_1$, 
	we obtain 
\[
	w \left( \theta \gamma(x) \right) \left( \gamma(y) \right) = 
	(bw_1 +aw_1 \circ w_0^{-1} \circ w_1) \left( \gamma(x) \right) \left( \gamma(y) \right).
\]
	Now the assertion of this lemma follows from the fact that $w_0, w_1 \in W \subset \Sym2 V^*$ 
	define symmetric bilinear forms on $V.$ \hfill \sq
\end{prf}
	Let $\mc{W} \subset \Sym2 L_W^*$ be the $k$-subspace consisting of $\widetilde{w} \in \Sym2 L_W^*$ satisfying (\ref{wdt}). 
	The dual space $L_W^*$ of $L_W$ is embedded into this space $\mc{W}$ by 
\begin{align*}
	D \colon L_W^* &\to \mc{W} \\
	\phi &\mapsto \left( (x, y) \mapsto \phi(xy) \right).
\end{align*} 
	The next lemma shows that all elements in $\mc{W}$ come from some elements of $L_W^*.$
\begin{lem}\label{Ne}
	\begin{enumerate}
	\item \label{cond1} The $k$-linear map $D \colon L_W^* \to \mc{W} $ is an isomorphism.
	\item Via the isomorphism of (\ref{cond1}), we consider $\mc{W}$ as an $L_W$-module.
		Then $\mc{W}$ is a free $L_W$-module of rank one, 
		and $\widetilde{w_0} = (\gamma^* \otimes \gamma^*)(w_0) \in \mc{W}$ generates $\mc{W}$ as an $L_W$-module.
	\end{enumerate}
\end{lem}
\begin{prf}
 	A $k$-bilinear form on $L_W$ can be regarded as an element of $\Hom_k \left( L_W \otimes_k L_W, k \right).$
 	And $k$-bilinear forms on $L_W$ satisfying (\ref{wdt}) are elements in 
 \[ 
 	\mc{W} = \Hom_k \left( L_W \otimes_{L_W} L_W, k \right) \isoc \Hom_k \left( L_W, k \right) = L_W^*.
 \]
 	The inverse of this isomorphism of $L_W$-modules is exactly $D.$ This proves the first assertion of Lemma \ref{Ne}.


	In order to prove the second assertion, we choose a $k^s$-basis of $L_W \otimes_k k^s.$
	Since $L_W$ is an \'etale $k$-algebra, there are $n$ different $k$-algebra homomorphisms 
	$\iota_i \colon L_W \rightarrow k^s (0 \leq i \leq n-1).$
	We put $\theta_i := \iota_i (\theta).$ By Proposition \ref{Reid},
	$\{\theta_i\}_{0 \leq i \leq n-1}$ are distinct, 
	and we can choose a $k^s$-basis $\{e_i\}$ of $L_W \otimes k^s$ such that for any $x_i \in k^s (0 \leq i \leq n-1),$
			\begin{align*}
				\widetilde{w_0} \left( \sum_{i=0}^{n-1} x_i e_i \right) &= \sum_{i=0}^{n-1} x_i^2, \\
				\widetilde{w_1} \left( \sum_{i=0}^{n-1} x_i e_i \right) &= \sum_{i=0}^{n-1} \theta_i x_i^2.
			\end{align*}
	Then we can check \[\theta(x_i e_i) = (w_0^{-1} \circ w_1)(x_i e_i) = \theta_i x_i e_i \in L_W \otimes_k k^s,\]
	and $L_W \otimes_k k^s$ coincides with the subspace of all elements in $\End (L_W \otimes k^s)$ having diagonal matrix representation
	with respect to the $k^s$-basis $\{e_i\}.$%

	Since $\{\theta_i\}_{0 \leq i \leq n-1}$ are distinct, the condition (\ref{wdt}) says all $q \in \mc{W} \otimes_k k^s$
	have diagonal matrix representation with respect to the basis $\{e_i\}.$ 
	This shows $\mc{W} \otimes k^s$ is a free $L_W \otimes_k k^s$-module of rank one generated by $\widetilde{w_0}.$
	Taking the Galois invariant subspaces, we prove the second assertion of Lemma \ref{Ne}.\hfill \sq
\end{prf}
	To summarize, we obtain the following diagram:
	$$\begin{xy}
	(0,0)*+\txt{$k$-basis} = "bs",
	(22,0)*+{\{\widetilde{w}_0,\widetilde{w}_1\}} = "wwbs",
	(47,0)*+{\{w_0,w_1\}} = "wbs",
	(72,0)*+{\{\id, \theta\}} = "mlbs",
	(0,10)*+\txt{2-dim. sp.} = "2d",
	(22,10)*+{\widetilde{W}} = "wtw",
	(47,10)*+{W} = "w",
	(72,10)*+{w_0^{-1}W} = "mcl",
	(0,20)*+\txt{$L_W$-module} = "lm",
	(22,20)*+{\mc{W} = \widetilde{w_0 L_W} } = "wwl",
	(47,20)*+{w_0 L_W} = "wl",
	(72,20)*+{L_W} = "l",
	(97,20)*+{V} = "v",
	(22,30)*+{\Sym2 L_W^*} = "sl",
	(47,30)*+{\Sym2 V^*} = "sv",
	(72,30)*+{\End V} = "ev", 
	(72,40)*+{\End L_W } = "el"
	\ar @{<->} "wwbs"; "wbs"
	\ar @{<->} "wbs"; "mlbs"
	\ar @{<->} "wtw"; "w"
	\ar @{<->} "w"; "mcl"
	\ar @{<->} "wwl"; "wl"
	\ar @{<->} "wl"; "l"
	\ar @{<->}^{\overset{\gamma}{\cong}} "l"; "v"
	\ar @{<->}^{\overset{\gamma^* \otimes \gamma^*}{\cong}} "sl"; "sv"
	\ar @{^{(}->}^{w_0^{-1}} "sv"; "ev"
	\ar @{}|{\cup} "wwl"; "wtw"
	\ar @{}|{\cup} "sl"; "wwl"
	\ar @{}|{\cup} "wtw"; "wwbs"
	\ar @{}|{\cup} "wl"; "w"
	\ar @{}|{\cup} "sv"; "wl"
	\ar @{}|{\cup} "w"; "wbs"
	\ar @{->}^{\cong}  "el";"ev"
	\ar @{}|{\cup} "ev"; "l"
	\ar @{}|{\cup} "l"; "mcl"
	\ar @{}|{\cup} "mcl"; "mlbs"
	\end{xy}$$
	
	As the final step of the construction of the map $\ms{orb}_n \rightarrow \ms{gc}_n,$
	we shall define $\alpha \in L_W^\times/k^\times L_W^{\times 2}.$ 
	In our condition, the $k$-algebra $L_W$ has the $k$-basis $\{\theta^i\}_{0 \leq i \leq n-1}.$ 
	Let $\{\theta^*_i\}_{0 \leq i \leq n-1}$ be the dual basis of $\{\theta^i\},$ and put 
	\[
		t_\theta:= D(\theta^*_{n-1}) \in \mc{W}.
	\]
	It is a nondegenerate quadratic form on $L_W$ because the matrix representation of $t_\theta$ with respect to the basis $\{\theta^i\}$ is 

\[ \begin{pmatrix}
0&0&\cdots&0&1\\
0&0&\cdots&1&*\\
\vdots&\vdots&\iddots&\vdots&\vdots\\
0&1&\cdots&*&*\\
1&*&\cdots&*&*
\end{pmatrix}_. \]

	By Lemma \ref{Ne}, for each $\widetilde{w} \in \mc{W},$ 
	there exists a unique element $\alpha_{w, \theta} \in L_W$ such that $\widetilde{w}= \alpha_{w, \theta} \cdot t_\theta.$
	Since $\widetilde{w_0}$ is nondegenerate, we obtain $\alpha \in L_W^\times$ satisfying $\widetilde{w_0} = \alpha \cdot t_\theta.$
	This $\alpha$ is what we need.

	In the construction of the map $\ms{orb}_n \rightarrow \ms{gc}_n$ as above, we have fixed two data:
	(i) a $k$-frame (i.e. ordered basis) $\{w_0, w_1\}$ of $W$, and
	(ii) an isomorphism of $L_W$-modules $\gamma \colon L_W \isoc V.$
	If we change $\gamma$ to $\gamma' = c \gamma$ for some $c \in L_W^\times,$ 
\begin{equation}\label{eqn}
\begin{split}
\left( (\gamma'^* \otimes \gamma'^*)(w) \right) (x,y) &= w \left( \gamma'(x), \gamma'(y) \right) \\
	&= w \left( c\gamma(x), c\gamma(y) \right) \\
	&= 	w \left( \gamma(x), c^2\gamma(y) \right) = \left( c^2 \cdot (\gamma^* \otimes \gamma^*)(w) \right) (x,y).
 \end{split}
\end{equation}
	Hence $\alpha$ is multiplied by an element of $L_W^{\times 2}.$
	To see the effects of a change of $k$-frame, we need another lemma.
	Recall again the trace form of an \'etale algebra is nondegenerate.
\begin{lem}[{\cite[III.6, Lemma 2]{Ser}}]\label{Sw}\mbox{}\\
	Let $\Tr_{L_W/k}, \N_{L_W/k}$ be the trace map and norm map from $L_W$ to $k$, 
	and \[P_\theta(t) := \N_{L_W/k}(t-\theta) \in k[t]\] the characteristic polynomial of $\theta.$ Then, we have
		\[\Tr_{L_W / k} \left( \frac{x}{P_\theta'(\theta)} \right) = \theta^*_{n-1}(x)\]
	for all $x \in L_W,$ where $P'_\theta(t) \in k[t]$ denotes the derivative of $P_\theta(t).$
\end{lem}
\begin{cor}\label{Sp}
	Let $\theta'$ be another generator of $L_W$ as a $k$-algebra.
	Then, we obtain $$t_{\theta'} =\frac{P_\theta'(\theta)}{P_{\theta'}'(\theta')} \cdot t_\theta.$$
\end{cor}

\begin{prf}\vspace{-5pt}
	\begin{equation*}
	\begin{split}
		\theta'^*_{n-1}(x) &= \Tr_{L_W / k} \left( \frac{x}{P_{\theta'}'(\theta')} \right)\\
		&= \Tr_{L_W / k} \left( \frac{P_{\theta}'(\theta)}{P_{\theta'}'(\theta')} \frac{x}{P_{\theta}'(\theta)} \right) 
		= \theta^*_{n-1}\left(  \frac{P_{\theta}'(\theta)}{P_{\theta'}'(\theta')} x \right). \mbox{\hfill \sq}
	\end{split}
	\end{equation*}
\end{prf}

	If we change the $k$-frame of $W$ to $\{dw_0+cw_1, bw_0+aw_1\}$ (and assume that $dw_0 + cw_1$ is nondegenerate),
	we see that $\theta$ changes to $\displaystyle \theta' := \frac{a \theta+b}{c \theta+d}$ because
	\[
	\begin{split}
		\theta' &= (dw_0 + cw_1)^{-1} \circ (bw_0 + aw_1) \\
		&= \left( w_0^{-1} \circ (dw_0 + cw_1) \right)^{-1} \circ w_0^{-1} \circ (bw_0 + aw_1) = (d+c\theta)^{-1}(b+a\theta).
	\end{split}
	\]
	Moreover, it is straightforward to see that
\begin{equation}\label{It}
	\frac{P_\theta'(\theta)}{P_{\theta'}'(\theta')} = \N_{L_W/k}(c\theta+d)\left(\frac{c\theta+d}{ad-bc}\right)^{n-1} \in k^\times L_W^{\times 2}
\end{equation}
	because $n$ is odd. By (\ref{eqn}), $\alpha$ is well-defined up to $k^\times L_W^{\times 2}.$
	Now the well-definedness of our map $\ms{orb}_n \rightarrow \ms{gc}_n$ is established.
\vspace{5pt}

	The construction of the inverse map $\ms{gc}_n \rightarrow \ms{orb}_n$ is easy;
	we send 
\[
	(L, \theta, \alpha) \mapsto (\alpha \cdot t_\theta, \alpha \theta \cdot t_\theta) \in k^2 \otimes \Sym2 L^*.
\]
	Via an isomorphism $L \isoc k^n$ as $k$-vector spaces, we obtain an element in $\ms{orb}_n.$
	Let us check the well-definedness. 
	A change $\alpha \mapsto a \alpha \, (a \in k^\times)$ does not affect the $k$-vector spaces $L$ and 
	the two-dimensional subspace $k \alpha \cdot t_\theta + k \alpha \theta \cdot t_\theta \subset \Sym2 L^*.$
	When one changes $\alpha \mapsto c^2\alpha \, (c \in L^\times),$ the isomorphism as $L$-module 
\begin{align*}
c \colon L &\isoc L \\ x &\mapsto cx
\end{align*}
	absorbs its effects.
	An argument similar to (\ref{It}) shows a change $\displaystyle \theta \mapsto \frac{a \theta+b}{c \theta+d}$ does not affect. 
	Hence this map is well-defined. Note that $(L, \theta, 1_L)$ is mapped to $(t_\theta, \theta \cdot t_\theta) \in \Sym2 L^*$ (see Section 2).
	
	We have to check these maps are inverses of each other. 
	First, we compare $W=(w_0, w_1)$ and $(\alpha \cdot t_\theta, \alpha \theta \cdot t_\theta)$ 
	obtained by the composition $\ms{orb}_n \rightarrow \ms{gc}_n \rightarrow \ms{orb}_n$.
	But an isomorphism as $L_W$-modules $\gamma: L_W \isoc V$ gives an isomorphism 
	between $(w_0, w_1)$ and $(\alpha \cdot t_\theta, \alpha \theta \cdot t_\theta).$
	
	Conversely, we take $(L, \theta, \alpha) \in \ms{gc}_n.$ 
	Because $(\alpha \cdot t_\theta)^{-1} \circ (\alpha \theta \cdot t_\theta) = \theta \in \End L ,$ we recover $L$ and the class of $\theta$. 
	Since $(t_\theta)^{-1} \circ (\alpha \cdot t_\theta) = \alpha \in \End V,$ 
	we recover $\alpha \in L^\times$ up to $k^\times L^{\times 2}.$
	
	In conclusion, we constructed the following bijection:
	
	\begin{equation}\label{second}
		\ms{orb}_n \overset{\sim}{\leftrightarrow} \ms{gc}_n = \left\{ \mbox{
		\begin{tabular}{l}
			equivalence classes of $(L, \theta, \alpha)$
		\end{tabular} } \right\}.
	\end{equation}
	
	Finally, we discuss the characteristic schemes. 
\begin{dfn}
	The {\it characteristic scheme} of a triple $(L, \theta, \alpha) \in \ms{gc}_n$ is 
	the projective equivalence class of the subscheme of $\p{1}_k$ 
	defined by the characteristic polynomial $P_\theta(t) = \N_{L/k}(t-\theta) \in k[t].$
\end{dfn}
	The projective equivalence between the characteristic schemes of $\ms{orb}_n$ and $\ms{gc}_n$ is given by $w_0.$
	It is easy to see that the bijection $\ms{orb}_n \isoc \ms{gc}_n$ preserves the characteristic schemes.
	
\subsection{Definition of $\ms{hec}_n$ and its relation to $\ms{gc}_n$}

\indent In this subsection, we define a set $\ms{hec}_n$ and construct a bijection 
between $\ms{gc}_n$ and $\ms{hec}_n$ preserving characteristic schemes.
To do this, we observe an isomorphism between the Jacobian of a hyperelliptic curve
 and a group scheme related to $\ms{gc}_n,$ and some data determining the isomorphism between them.

Recall that $m \geq 1$ is a positive integer and $n=2m+1.$
Let $P$ be a $k$-rational point of $\p1_k$. Take a nonsingular projective hyperelliptic curve $C$ over $k$ of genus $m$ 
with a hyperelliptic involution $\iota_C \colon C \rightarrow C$. The involution induces a morphism $\pi \colon C \to \p1_k.$
We assume there exists a $k$-rational point $\widetilde{P} \in C(k)$ 
such that $\iota_C(\widetilde{P})=\widetilde{P}$, and $\pi(\widetilde{P})=P$ (i.e. $\widetilde{P}$ is a Weierstrass point of $C$).
This is equivalent to say $\pi$ is ramified at $P \in \p1_k(k).$ 

Let $T_C$ be the an $(n+1)$-points $k$-subscheme of $\p1_k$ where the quotient morphism $\pi$ is ramified.
This scheme $T_C \subset \p1_k$ is a disjoint union of the $k$-rational point $P$ and an $n$-points subscheme $S_C$ (i.e. $T=\{P\} \sqcup S_C$).
Then the coordinate ring $L_C:=H^0 \left( S_C, \Oh_{S_C}\right)$ is an \'etale $k$-algebra of degree $n$.

When we take a basis $\{ 1, t \}$ of $H^0 \big( C, \Oh(2\widetilde{P}) \big) = H^0 \big( \p1_k, \Oh(P) \big)$,
 we obtain a hyperelliptic equation $y^2=f(t)$ of $C$.
This $f(t)$ is a polynomial in $k[t]$ of degree $n$ such that $L_C \cong k[t]/(f(t)).$ On the other hand, by the inclusion $S_C \hookrightarrow \p1_k$, 
we have a $k$-module homomorphism 
\begin{equation}\label{morph}
	H^0\left( \p1_k, \Oh(P) \right) \to H^0\left( S_C, \Oh_{S_C}\right) = L_C.
\end{equation}
If we write the image of $t$ as $\theta,$ the characteristic polynomial of $\theta$ coincides with $f(t).$

On this step, we have the following isomorphism:

\begin{prop}[c.f. {\cite[Theorem1.1]{Scha}}]\label{Fi}
	Let $J_C$ be the Jacobian variety of $C$, and $J_C[2]$ the group scheme of 2-torsion points on $J_C.$
	Then there exists an isomorphism of group schemes over $k$ : which depends on $C$ and its hyperelliptic equation:
	\[ w_C: J_C[2] \isoc \Res_{L_C/k}(\mu_2)/\Delta(\mu_2).\] 
\end{prop}

	From the long exact sequence of Galois cohomology, we obtain isomorphisms; 
	\[
		H^1\left(G_k, J_C[2] \right) \cong H^1\left( G_k, \Res_{L_C/k}(\mu_2)/\Delta(\mu_2) \right)
		\cong L_C^\times / k^\times L_C^{\times 2}, 
	\]
	where the second isomorphism is given by Kummer theory, and the fact that $n$ is odd.
	
	Now we define the set $\ms{hec}_n$ related to hyperelliptic curves. 
\begin{dfn}\label{hec}
	Let us consider the following triple:
	\begin{itemize}
		\setlength{\parskip}{0pt}
		\setlength{\itemsep}{0.0pt}
		\item $C$ is a nonsingular projective hyperelliptic curve over $k$ of genus $m$,
		\item $P \in T_C(k)$ is the image of a $k$-rational Weierstrass point of $C$,
		\item $\eta$ is an element of $H^1\left( G_k, J_C[2] \right)$.
	\end{itemize}
	We say two such triples $(C, P, \eta), (C', P', \eta')$ are {\it equivalent} if the following conditions hold:\\ 
 (i) there is an automorphism $\Phi \colon \p1_k \isoc \p1_k$ which induces an isomorphism $\Phi |_{S_{C}} \colon S_{C} \isoc S_{C'},$ \\
 (ii) $\eta$ is mapped to $\eta'$ through the following isomorphisms (see Proposition \ref{Fi}) : 
\vspace{-3pt} 
\[ H^1\left( G_k, J_C[2] \right) \cong L_C^\times / k^\times L_C^{\times 2} \overset{\phi}{\overset{\sim}{\leftarrow}}
 L_{C'}^\times / k^\times L_{C'}^{\times 2} \cong H^1\left( G_k, J_{C'}[2] \right)\]
 where $\phi$ is induced by $\Phi |_{S_{C}}.$ 
 
 We define $\ms{hec}_n$ as the set of equivalence classes of triples as above. 
The {\it characteristic scheme} of the triple $(C, P, \eta)$ is the projective equivalence class of $S_C \subset \p{1}_k.$
\end{dfn}

	Now we shall construct a bijection between $\ms{gc}_n$ and $\ms{hec}_n$ preserving characteristic schemes.
	
	First, we take a representative $(C, P, \eta) \in \ms{hec}_n.$ 
	The way to get the corresponding element $(L, \theta, \alpha) \in \ms{gc}_n$ is already explained in the beginning of this subsection.
	The algebra $L$ is $H^0\left( S_C, \Oh_{S_C} \right).$ 
	Fixing a basis $\{1, t \}$ of $H^0\big( \p1_k,\Oh_{\p1_k}(P)\big),$ 
	we define $\theta$ as the image of $t$ by the restriction map 
	$H^0\big( \p1_k,\Oh_{\p1_k}(P)\big) \to H^0\left( S_C, \Oh_{S_C} \right).$ 
	And we take a hyperelliptic equation $y^2 = f(t)$ of $C$, which gives the isomorphism $w_C$ in Proposition \ref{Fi}.
	Then we obtain $\alpha \in L^\times / k^\times L^{\times 2}$ as the image of $\eta \in H^1\left( G_k, J_C[2] \right).$
	
	We must check the well-definedness of this map. First, we define an isomorphism (rather than equivalence)
	between triples $(C, P, \eta), (C', P', \eta')$ as a $k$-isomorphism $C \isoc C'$ whose inducing morphisms send $P$ to $P'$, and $\eta$ to $\eta'.$
	We must show the above map is well-defined as a map from the set of equivalent classes of triples to $\ms{gc}_n.$
	It is enough to consider the following two cases:
	\begin{itemize}\label{cond}
     \setlength{\parskip}{0pt}
     \setlength{\itemsep}{0.0pt}
     \item Fix $C$ and $P$, but change the basis $\{1, t \}$ to $\{1, at+b\} \; (a \in k^\times, b \in k).$
     \item Change the polynomial $f(t)$ to $a^2 f(t)$ for some $a \in k^{\times}$.
    \end{itemize}
    All of these changes do not affect $L.$ Changes of second kind also do not affect $\theta,$
    while those of first kind affect $\theta$ as M\"obius transformations.
    So it suffices to show that these changes do not affect $\alpha.$
    
    To do this, we recall the definition of $w_C$ in \cite{Scha}. Put $L^s := L \otimes_k k^s$.
     We have \[ L^s \isoc (k^s)^n \quad ; \quad \theta \mapsto \left( \theta_1, \theta_2, \cdots, \theta_n \right). \]
     We put \[ \overline{P} := \left( (\theta_1, 0) - P, (\theta_2, 0) - P, \cdots, (\theta_n, 0) - P) \right) \in \Res_{L/k}\left( J_C[2] \right)(k^s). \]
     For any $R \in J_C[2](k^s),$ we define $w_C(R)$ explicitly as
\[
     w_C(R):= e_2 \left( \Delta(R), \overline{P} \right) = 
	\left( e_2\left( R, (\theta_i, 0)-P \right)_i \right) \in \left( \Res_{L/k}(\mu_2) / \Delta(\mu_2) \right) (k^s)
\]
     where $\Delta$ is the diagonal embedding, and $e_2$ is the Weil pairing for 2-torsion subgroup.
     From this description, we see that the two changes do not affect $w_C,$ and therefore $\alpha.$
     
     Now what remains to show is as follows. 
	Take two equivalent non-isomorphic triples $(C, P, \eta)$ and $(C', P', \eta'),$
	and assume $S_C = S_C'$ in $\p1_k.$ Do these triples correspond to the same element in $\ms{gc}_n$?

	Since $S_C = S_C',$ we identify $L_C =L_{C'}$ and $\alpha = \alpha'.$ 
	Only the possible difference is $\theta$ and $\theta'.$ It is easy to see that $\theta'$ is a M\"obius transformation of $\theta.$
	Hence $\ms{hec}_n \to \ms{gc}_n$ is well-defined.
	
\vspace{5pt}
	Second, we take $(L, \theta, \alpha) \in \ms{gc}_n.$ Let $P_\theta(t)$ be the characteristic polynomial of $\theta.$
	We put $C \colon y^2=P_\theta(t), P = \{ \infty \},$ and $\eta$ corresponds to $\alpha$ via $w_C.$
	Again we must check the well-definedness of this map.
	
	Recall the equivalence relation defining $\ms{gc}_n.$ It is a combination of transformations of two kinds.
	First kind is an isomorphism $L \isoc L'$ which maps $\theta$ to $\theta'$, $\alpha$ to $\alpha'.$
	Second kind is a M\"obius transformation $\theta \mapsto \displaystyle \frac{a\theta+b}{c\theta+d}.$
	A transformation of first kind induces an isomorphism $C \isoc C'$ which maps $P$ to $P',$ and  $\eta$ to $\eta'.$ 
	So we may consider only a transformation of second kind.
	But in this case, an automorphism 
	\begin{align*}
		\phi \colon \p1_k &\isoc \p1_k \\
		(x:z) &\mapsto (-bz+dx: az-cx)
	\end{align*}
	induces projective equivalence $S_C \isoc S_{C'},$ which sends $P'$ to $(d: -c)$ and $\alpha$ to $\alpha'.$
	So we may identify $S_C = S_{C'}, P'=(d: -c), \alpha=\alpha'.$ This shows the well-definedness of a map $\ms{gc}_n \to \ms{hec}_n.$
	It is straightforward to see that these two maps are inverses of each other.
	

	The required property that the maps preserve characteristic schemes is clear from our construction. In conclusion, we have constructed a bijection:
\begin{equation}\label{third}
	\ms{gc}_n = \left\{ \hspace{-5pt} \mbox{
	\begin{tabular}{l}
		equivalence classes of $(L, \theta, \alpha)$
	\end{tabular} } \hspace{-5pt} \right\} \overset{\sim}{\leftrightarrow} \left\{ \hspace{-5pt} \mbox{
	\begin{tabular}{l}
		equivalence classes of $(C, P, \eta)$
	\end{tabular} } \hspace{-5pt} \right\} = \ms{hec}_n
\end{equation}
preserving characteristic schemes.

\begin{exa}\label{exa}
	(i) For a fixed $(C, P),$ the map $ H^1(G_k, J_C[2]) \to \ms{hec}_n $ is not injective in general.
	Let us show the following example.

	Take $C \colon y^2=t^5-1$ over $k = \Q (\zeta_5),$ 
	where $\zeta_5$ is a primitive fifth root of unity in $\overline{\Q}.$
	The $k$-automorphism $\Phi \colon (s : t) \mapsto (\zeta_5s : t)$ of $\p1_k$ induces automorphisms 
	$S_C \isoc S_C,$ and 
		\begin{align*}
			\phi \colon L_C^\times / k^\times L_C^{\times 2} &\isoc L_C^\times / k^\times L_C^{\times 2} \\
			\sum_{i=0}^{n-1} a_i \theta^i &\mapsto \sum_{i=0}^{n-1} a_i \zeta_5^i \theta^i. 
		\end{align*}
	And we can find $\alpha \in L_C^\times = \left( k[t]/(f(t)) \right) ^\times$
	such that $\alpha \neq \phi(\alpha)$ in $L^\times / k^\times L^{\times 2}.$ 
	Since two triples $(L_C, \theta, \alpha)$ and $\left( L_C, \theta, \phi(\alpha) \right)$ are equivalent, 
	they define the same element in $\ms{hec}_5.$
	We can determine the fiber of this map in Section 2.\\	
	(ii) Note that two triples $(C, P, \eta), (C', P', \eta')$ can be equivalent in the sense of Definition \ref{hec}
	even if $C, C'$ are not $\overline{k}$-isomorphic (e.g. $y^2 = t^5-1$ and $y^2=(t+1)^5+1$ over $k=\Q$ or $\Q(\zeta_5)$). 
	
\end{exa}

\section{Definition of the map $i_C$ and the proof of the main theorem}
 Now we shall define the map $i_C$ in our main theorem (Theorem \ref{main}). 
	Let $f(t) \in k[t]$ be a separable polynomial of degree $n$, and $C$ a nonsingular hyperelliptic curve over $k$ which has an affine model $y^2=f(t)$. 
	Let $P = \infty \in \p1_k$ be the point at infinity. 
	Put $S_C :=  \Spec k[t]/(f(t))$. The double covering $\pi \colon C \rightarrow \p1_k$ is ramified at $\{ \infty \} \cup S_C \subset \p1_k.$
	Then we define 
	\[
	i_C \colon H^1\left( G_k, J_C[2] \right) \to \ms{ciq}_n
	\]
	as the composite
	\begin{alignat*}{8}
		H^1(G_k, J_C[2]) \; &\to \; \ms{hec}_n& \; &\to \;  \ms{gc}_n& \; &\to \;  \ms{orb}_n &\; &\to \;  \ms{ciq}_n &\\
		\eta \; &\mapsto \; (C, P, \eta)& \; &\mapsto \; (L, \theta, \alpha)& \; &\mapsto \; (A,B)& \; &\mapsto \; [X] &=: i_C(\eta).
	\end{alignat*}
	where $[X]$ denotes the projective equivalence class of nonsingular polarized complete intersection of two quadrics over $k$ corresponding to $(C, P, \eta).$
	Recall that the isomorphism in Proposition \ref{Fi} depends on the choice of $P,$ or $\theta.$
	The following proposition determine the image of $i_C.$
 
 \begin{prop}\label{main1}
	The image of $i_C$ consists of all elements $[X] \in \ms{ciq}_n$ whose characteristic scheme is projectively equivalent to $S_C.$
 \end{prop}
 \begin{prf}
	This follows from the results in the previous section because our bijections preserve the characteristic schemes. \hfill \sq
 \end{prf}

	The map $i_C$ in Theorem \ref{main} is not injective in general (see Example \ref{exa} (i)).
	We say $\theta' \in L$ are {\it conjugate} to $\theta$ over $k$ if they have the same characteristic polynomial over $k.$
	The following proposition determines the fiber of $i_C.$

\begin{prop}\label{Ru}
	Take $[X] \in \ms{ciq}_n,$ and let $\displaystyle (L, \theta, \alpha)$ be the triple corresponding to $[X].$ 
	Let us write $\alpha$ as
	\[
		\alpha = \sum_{i=0}^{n-1} a_i \theta^i.
	\] 
	We put 
	\[
		\Aut_{L, \theta}:= \left\{A \in \GL_2(k) \; \mid \; A \cdot \theta \; \mbox{{\rm is conjugate to }} \theta \; \mbox{{\rm over }} \; k \right\}.
	\]
	Then the fiber $i_C^{-1}([X])$ is bijective to 
	\begin{equation}\label{orbit}
		\left\{ \left. \sum_{i=0}^{n-1} a_i (A \cdot \theta)^i \in L^\times / k^\times L^{\times 2} \right|  A \in \Aut_{L, \theta} \right\}_.
	\end{equation}
\end{prop}

\begin{prf}
	We can define a map $j_C: L^\times / k^\times L^{\times 2} \to \ms{gc}_n$ which makes the following diagram commute:
\[ \begin{CD}
	H^1(G_k, J_C[2]) @>{i_C}>> \ms{hec}_n\\
	@V{w_C}V{\cong}V @VV{\cong}V\\
	L^\times / k^\times L^{\times 2} @>{j_C}>> \ms{gc}_n.
\end{CD} \] 
	Since the vertical arrows are bijective, we may consider the map $j_C.$
	
	The fiber of $j_C$ is 
	\[
		\{ \beta \in L^\times / k^\times L^{\times 2} \mid (L, \theta, \alpha) \mbox{ is equivalent to } (L, \theta, \beta) \}.
	\]
	If $(L, \theta, \alpha)$ is equivalent to $(L, \theta, \beta),$ there exists a $k$-algebra automorphism $\phi \colon L \isoc L$
	which maps $\alpha$ to $\beta$ and $\theta$ to a M\"obius transformation $A \cdot \theta$ for some $A \in \GL_2(k).$
	This shows $\beta$ is described as in (\ref{orbit}) for $A \in \Aut_{L, \theta}.$

	Conversely, take $A \in \Aut_{L, \theta}$ and a triple $(L, \theta, \alpha).$ 
	Put 
	\[
		\beta = \sum_{i=0}^{n-1} a_i (A \cdot \theta)^i .
	\]
	Then by an automorphism of $k$-algebra $L \isoc L$ satisfying $\theta \mapsto A \cdot \theta,$ the triple $(L, \theta, \beta) $ is equivalent to $(L, \theta, \alpha).$  \hfill \sq
\end{prf}

	In the rest of this section, we will study elements in each sets corresponding to the identity elements in Galois cohomology.
	Such elements are called quasi-split (see Remark \ref{rm}).
\begin{dfn}
	A triple $(C, P, \eta) \in \ms{hec}_n$ is called {\it quasi-split} \. if $\eta =0 \in H^1\left(G_k, J_C[2] \right).$ 
	An element of $\ms{gc}_n, \ms{orb}_n, \ms{ciq}_n$ is called {\it quasi-split} if the corresponding element in $\ms{hec}_n$ is quasi-split.
\end{dfn}
\begin{rmk}\label{rm}
	In \cite{BG}, quasi-split objects in $\ms{orb}_n$ are called {\it distinguished orbits}. 
	And in \cite{Skor2}, quasi-split objects in $\ms{gc}_5$ are called {\it quasi-split del Pezzo surfaces} (see next section).
\end{rmk}
	A triple $(L, \theta, \alpha) \in \ms{gc}_n$ is quasi-split if and only if $\alpha = 1_L$ in $L^\times / k^\times L^{\times 2}$ (see Proposition \ref{Fi}).
	From this, we obtain the following corollary.
\begin{cor}
	Let $i_C \colon H^1(G_k, J_C[2]) \to \ms{ciq}_n$ be the map defined in the beginning of this section.
	For a quasi-split element $[X] \in \ms{ciq}_n,$ we have $i_C^{-1}([X])={0}.$
\end{cor}
	What is a characterization of quasi-split objects in $\ms{orb}_n, \ms{ciq}_n?$ 
	By the bijection we constructed, $(L, \theta, 1_L) \in \ms{gc}_n$ maps to $(L, kt_\theta+k\theta \cdot t_\theta) \in \ms{orb}_n.$ 
	By definition of $t_\theta$, the pencil of quadratic forms corresponding to $(L, kt_\theta + k\theta \cdot t_\theta)$
	has a common isotropic $k$-subspace of dimension $m,$ 
\[
	\mc{L}^{m-1} := \left\{ a_0 + a_1 \theta + \dots + a_{m-1} \theta^{m-1} \mid a_i \in k \quad (0 \leq i \leq m-1) \right\},
\]
	because for any $u, v \in \mc{L}^{m-1}$ and  $a, b \in k,$ $(a+b\theta)uv$ is a polynomial in $\theta$ of degree $\leq 2m-1 = n-2,$
	and hence 
	\[
		\theta^*_{n-1}((a+b\theta)uv)=0.
	\]
	Indeed, this property characterizes the quasi-split objects in $\ms{orb}_n.$ 
\begin{prop}
	Let $w = \alpha \cdot t_\theta \in \Sym2 L^*$ be a nondegenerate quadratic form on $L.$ 
	Assume that there exists an $m$-dimensional $k$-subspace $M$ of $L$ which is isotropic for $w$ and $\theta \cdot w.$
	Then $\alpha = 1_L$ in $L^\times / k^\times L^{\times 2}$.
\end{prop}
\begin{prf}
	First, we prove the following lemma:
	\begin{lem}\label{crit}
		If $\ell \in M$ satisfies $\ell \neq 0$ and $\ell \theta^i \in M (0 \leq i \leq m-1),$
		the elements $\{ \ell \theta^i \}_{0 \leq i \leq m-1}$ form a $k$-basis of $M$.
		Moreover, $\ell \in L^\times$ and $\alpha = 1_L$ in $L^\times / k^\times L^{\times 2}$.
	\end{lem}
	\begin{prf}
		By assumption, $\ell^{2} \cdot w$ and $\ell^{2} \theta \cdot w$ have 
		isotropic $k$-subspace $\left< 1, \theta, \cdots, \theta^{m-1} \right>$ because
		\begin{align*}
			(\ell^2 \cdot w)(\theta^i, \theta^j) &= w(\ell \theta^i, \ell \theta^j) = 0 \\
			(\ell^2 \theta \cdot w)(\theta^i, \theta^j) &= w(\ell \theta^{i+1}, \ell \theta^j) = 0
		\end{align*}
		for $0 \leq i, j \leq m-1.$
		This shows the $k$-linear form $w(\cdot, \ell^2)$ kills $1, \theta, \cdots, \theta^{2m-1} = \theta^{n-2}.$
		Hence $w(x, \ell^2) = a \theta^*_{n-1} (x),$ so \[ \ell^2 \cdot w = \alpha \ell^2 \cdot t_\theta = a \cdot t_\theta\] for some $a \in k.$
		By Lemma \ref{Ne}, we deduce $\alpha \ell^2 = a.$

		Since $L$ is an \'etale $k$-algebra and $\ell \neq 0,$ we have $\ell^2 \neq 0.$ 
		And since $w$ is nondegenerate, we also obtain $\alpha \in L^\times.$ Hence $\alpha \ell^2 = a \neq 0.$
		This shows $\alpha = 1_L$ in $L^\times / k^\times L^{\times 2}$
		and $\ell \in L^\times.$ Since $\{ \theta^i \}_{0 \leq i \leq n-1}$ is linearly independent over $k$ and $\ell \in L^\times,$
		$\{ \ell \theta^i \}_{0 \leq i \leq n-1}$ is linearly independent over $k.$
		By comparing dimensions, \[ M = k\ell + k\ell \theta + \cdots + k\ell \theta^{m-1}, \]
		and the elements $\{ \ell \theta^i \}_{0 \leq i \leq m-1}$ form a basis of $M$. \hfill \sq
	\end{prf}

	By this lemma, we have only to show there exists an element $\ell \in M$ satisfying the assumption of Lemma \ref{crit}.
	
	Before discussion, we prepare some notation.
	For a $k$-subspace $N \subset L$, we write the orthogonal complement of $N$ with respect to $w$ by
\[
N^{\bot} := \left\{ a \in L \mid w(a, \ell) = 0 \; \mbox{ for all } \ell \in N \right\}.
\] 
	Also we put $a N := \left\{ a n \mid n \in N \right\}$ for $a \in L.$
	By assumption, $M^\bot$ contains $M$ and $\theta M.$ 
	
	Next consider the $k$-linear map $\theta \times : M \to \theta M$ defined by $x \mapsto \theta x.$
	Assume $\ker (\theta \times) \neq 0$ and take a non-zero element $\ell' \in \ker(\theta \times).$
	Then, $\ell'$ satisfies the condition of $\ell$ in Lemma \ref{crit}.
	But $\ell' \not\in L^\times$ because $\theta \ell' =0,$ and we get a contradiction. 
	Therefore, $\ker (\theta \times)=0$ and $M$ is isomorphic to $\theta M$ via the map $\theta \times.$
	
	Put $M_j := M \cap \theta M \cap \cdots \cap \theta^j M$ and $M_0 =M.$
	We shall show that $\dim_k M_j \geq m-j$ inductively.
	If it is proved, $\dim_k M_{m-1} \geq 1,$ and a non-zero element $\ell \in M_{m-1}$ satisfies the assumption of Lemma \ref{crit}.
	
	As a first step, we have $\dim_k M^\bot = m+1$ because $w$ is nondegenerate and $n=2m+1.$
	Since $\dim_k M = \dim_k \theta M = m,$ we conclude 
	\begin{align*}
		\dim_k M_1 &= \dim_k \left( M \cap \theta  M \right) \\
		&= \dim_k M + \dim_k \theta M - \dim_k (M+\theta M)\\
		&\geq \dim_k M + \dim_k \theta M - \dim_k M^\bot\\
		&= m-1.
	\end{align*}
	If $\dim_k M_1 = m,$ then $M = \theta M,$ and $\dim_k M_j = m$ for all $j \geq 0.$
	So we may assume $\dim_k M_1 = m-1.$
	
	By induction, we may assume $\dim_k M_i = m-i \; (i \leq j) $ for some $j \geq 1.$ The space $M_{j-1}$ contains $M_j$ and $\theta M_j,$
	and we observe $M_j \cap \theta M_j = M_{j+1}.$ As above, we conclude 
	\begin{align*}
		\dim_k M_{j+1} 
		&= \dim_k M_j + \dim_k \theta M_j - \dim_k (M_{j}+\theta M_j)\\
		&\geq  \dim_k M_j + \dim_k \theta M_j - \dim_k M_{j-1}\\
		&= m-j-1.
	\end{align*}
	If $\dim_k M_{j+1}=m-j,$ then $M_j = \theta M_j,$ and $\dim_k M_i = m-j$ for all $i \geq j.$
	So we may assume $\dim_k M_j = m-j-1,$ and the induction proceeds. \hfill \sq

\end{prf}
	To summarize, we have the following theorem.
\begin{thm}\label{quasi}
	An element in $\ms{orb}_n$ is quasi-split if and only if the pencil of quadratic forms $W \subset \Sym2 V^*$ corresponding to it
	has a common isotropic $m$-dimensional $k$-subspace of dimension $m.$
	An element $[X] \in \ms{ciq}_n$ is quasi-split if and only if $X$ contains a linear $k$-subvariety isomorphic to $\p{m-1}_k$.
\end{thm}
\begin{cor}\label{quasi2}
	We fix a reduced $n$ points subscheme $S$ of $\p1_k$. 
	There is a unique projective equivalence class of nonsingular complete intersection of two quadrics $X$ in $\p{n-1}_k$ satisfying the following properties.
	\begin{itemize}
		\setlength{\parskip}{0pt}
		\setlength{\itemsep}{0.0pt}
		\item The characteristic scheme $X$ is projectively equivalent to $S.$
		\item The variety $X$ contains a linear $k$-subvariety isomorphic to $\p{m-1}_k$.
	\end{itemize}
\end{cor} 
	Theorem \ref{quasi} and Corollary \ref{quasi2} finish the proof of Theorem \ref{main}
	because quasi-split objects in $\ms{ciq}_n$ can be written as $i_C(0)$ for some $C.$

	Finally, we shall show Theorem \ref{main2}. 

\begin{thm}\label{Jp}
	Let $\delta \colon J_C(k) \to H^1(G_k, J_C[2])$ be the connecting homomorphism in Galois cohomology.
	Let $[X] \in \ms{ciq}_n$ be an element contained in the image of the composite map $i_C \circ \delta$ for some $C.$
	Then the variety $X$ contains a linear $k$-subvariety isomorphic to $\p{\lfloor \frac{m-1}{2} \rfloor}.$ 
\end{thm}
	\begin{prf}
	In the following, we use the results ion the Jacobian variety of a hyperelliptic curve summarized in Bruin--Flynn's paper (\cite[Section5]{BF}).
	
	Take an element $D \in J_C(k).$
	The divisor class $D$ is represented by a divisor of the following form  
\[
	\{g(x) = 0, y= h(x) \} \cdot C - \left( \mathrm{deg}(g) \right)\infty.
\]
	Here $g(x), h(x) \in k[x], $ and $0 \leq \mathrm{deg} (h) < \mathrm{deg} (g) \leq m.$
	
	Then via the composite $w_C \circ \delta$, where $w_C$ is the map in Proposition \ref{Fi},
	the element $D \in J(k)$ is mapped to 
\[
	\alpha_D := (-1)^{\mathrm{deg}(g)}g(\theta) \in L^\times / k^\times L^{\times 2}.
\]
	(In fact, in \cite{BF}, Bruin and Flynn use another group $\Ker \left( L^\times / L^{\times 2} \to k^\times / k^{\times 2} \right).$
	But it is easy to justify our argument using the isomorphism
	\begin{align*}
		\Ker (\N_{L/k} \colon L^\times / L^{\times 2} \to k^\times / k^{\times 2} ) &\isoc L^\times / k^\times L^{\times 2} \\
		\alpha &\mapsto \alpha.
	\end{align*} Here we use the assumption that $[L : k] = n$ is odd.)
	Clearly, 
	\[
		M = k + k\theta + \cdots + k\theta^{\lfloor \frac{m-1}{2} \rfloor}. 
	\]
	 is a common isotropic $k$-subspace of dimension $\lfloor \frac{m-1}{2} \rfloor+1$ 
	for $\alpha_D \cdot t_\theta$ and $\alpha_D \theta \cdot t_\theta.$ \hfill \sq
	\end{prf}
	\begin{rmk}
		Theorem \ref{Jp} only gives a necessary condition for an element $[X] \in \ms{ciq}_n$ to be in the image of $i_C \circ \delta.$
		This is not a sufficient condition.
	\end{rmk}

\section{Examples}
	\subsection{Case for $n=3$: Four points subschemes of $\p2_k$}
	This case is well-known. The nonsingular complete intersections of two quadrics in $\p2_k$ are four points subschemes of $\p2_k$ in general position (i.e.\ no three points are geometrically collinear). From Proposition \ref{Reid}, we can check the converse is true.
	In fact, if one takes a four points $k$-subscheme of $\p2_{k}$ in general position, we can regard four points are $[1: \pm1: \pm1]$ in some coordinate system over $k^s$. So quadrics through such four points consist a two-dimensional subspace of $\Sym2 k^3$.
	This shows four points subschemes of $\p2_k$ in general position are nonsingular complete intersections of two quadrics in $\p2_k.$
	We say two four points subschemes $X, Y \subset \p2_k$ are projectively equivalent if there exists an automorphism of $\p2_k$ over $k$
	sending $X$ to $Y.$ In conclusion, we obtain:
	\begin{prop}
		There is a bijection between projective equivalence classes of four points $k$-subschemes of $\p2_k$
		in general position and projective equivalence classes of nonsingular pencils of quadratic forms in $\Sym2 k^3$.
	\end{prop}

	A four points subscheme of $\p{2}_k$ is {\it quasi-split} if and only if it contains a $k$-rational point.

	\subsection{Case for $n=5$: del Pezzo surfaces of degree four}
	We recall the definition of del Pezzo surfaces and some properties of them. For more details, see \cite{Dolg} for example.
	
	A {\it del Pezzo surface} over $k$ is a nonsingular projective surface $X$ over $k$ with ample anticanonical divisor $-K_X$. Its {\it degree} is the self-intersection number $(-K_X)^2$. From now on, we concentrate on the case of degree four.
	
	Let $X$ be a del Pezzo surface of degree four over $k$. We know that $V = H^0\left( X, -K_X \right)$
	is a five-dimensional $k$-vector space, and $X$ is realized as a nonsingular complete intersection of
	two quadrics in $\p{}(V) \cong \p4_k.$ Conversely, all pairs of two quadrics in $\p4_k$ 
	which intersect completely and smoothly define del Pezzo surfaces of degree four.
	
	An isomorphism of del Pezzo surfaces $X \isoc Y$ induces an isomorphism of $k$-vector spaces 
	$H^0 \left( X, -K_X \right) \isoc H^0\left( Y, -K_Y \right)$. 
	Conversely, since $H^0\left( X, -K_X \right)$ is a complete linear system, 
	all isomorphisms of del Pezzo surfaces of degree four over $k$ are induced by such $k$-linear isomorphisms. 
	
	Furthermore, the kernel of the map 
\[
\cdot |_X \colon  \Sym2{} \left( H^0 \left( X, -K_X \right) \right)  \to H^0 \left( X, -2K_X \right)
\]
 is a two-dimensional $k$-vector space \cite{Dolg}. 
	Hence the set of isomorphism classes of del Pezzo surfaces of degree four 
	can be identified with the set of 
	projective equivalence classes of pencils of quadrics $\Phi$ in $\p4_k.$ 
	This shows $\ms{ciq}_5$ coincides with the set of isomorphism classes of del Pezzo surfaces of degree four over $k$. So we conclude:
	\begin{prop}\label{dP4}
		There is a bijection between isomorphism classes of del Pezzo surfaces of degree four over $k$
	and projective equivalence classes of nonsingular pencils of quadratic forms in $\Sym2 k^5$.
	\end{prop}
	
	A del Pezzo surface of degree four $X$ is called {\it quasi-split} when it contains a $k$-line. 
	In \cite{Skor2}, Skorobogatov showed that a del Pezzo surface of degree four $X$ defined over $k$
	is quasi-split if and only if $X$ is isomorphic to the blow-up of $\p2_k$ at a five points subscheme defined over $k$.
	
	In  \cite[Section 2]{Skor2}, Skorobogatov gives a bijection between isomorphism classes of del Pezzo surfaces of degree four over $k$
	and another set: pairs of projective equivalence classes of five points schemes $S \subset \p1_k$ 
	and elements $\lambda \in H^1(G_k, \Res_{L/k}\mu_2/\Delta(\mu_2))$ where $L=H^0(S, \Oh_S).$ 
	We write this set as $\ms{gc'}_5$.
	In his bijection, a pair $(S, \lambda)$ corresponds to
	\[
		X := \left\{ x \in \p{}(\Res_{L/k}\mathbb{A}^1_L) (\cong \p4_k) \; \left| \; \Tr_{L / k} \left( \frac{\lambda x^2}{P_\theta'(\theta)} \right) 
		= \Tr_{L / k} \left( \frac{\lambda \theta x^2}{P_\theta'(\theta)} \right) =0 \right. \right\}
	\]
	where $\theta \in L$ is the image of $t \in H^0(\p1_k, \Oh(\infty))$ via the homomorphism in (\ref{morph}).
	
	As we have seen, the isomorphism class of this variety $X$ corresponds to $(L, \theta, \lambda) \in \ms{gc}_5.$
	Hence we have a bijection between $\ms{gc}_5$ and $\ms{gc'}_5.$ 
	If $(L, \theta, \lambda) \in \ms{gc}_5$ corresponds to $(S, \lambda) \in \ms{gc'}_5,$
	$S$ is the characteristic scheme of $(L, \theta, \lambda).$ 
	Our bijection is compatible with Skorobogatov's bijection in this way. 
	In particular, quasi-split del Pezzo surfaces of degree four in \cite{Skor2} corresponds to quasi-split objects.
	

{\sc Graduate School of Mathematics, Faculty of Science, Kyoto University,
Kyoto 606-8502, Japan}

{\it E-mail address}: {\tt yasu-ishi@math.kyoto-u.ac.jp}

\end{document}